# Power sum identities with generalized Stirling numbers

Khristo N. Boyadzhiev
Department of Mathematics
Ohio Northern University
Ada, OH 45810, USA
k-boyadzhiev@onu.edu

**Abstract** We prove several combinatorial identities involving Stirling functions of the second kind with a complex variable. The identities involve also Stirling numbers of the first kind, binomial coefficients and harmonic numbers.

## 1. Introduction

Butzer, Kilbas and Trujillo [2] defined the Stirling functions of the second kind by

$$S(\alpha,k) = \frac{1}{k!}\sum_{j=1}^{k} (-1)^{k-j}\binom{k}{j} j^\alpha, \qquad (1)$$

for all complex numbers $\alpha \neq 0$ and all positive integers $k$. This definition is consistent with the definition given by Flajolet and Prodinger [5]. When $\alpha = n$ is a positive integer, $S(n,k)$ are the classical Stirling numbers of the second kind [2], [3], [8]. The purpose of this note is to prove the five power sum identities (6), (17), (20), (23) and (24) below involving the Stirling functions $S(\alpha,k)$. In fact, we describe a general method for obtaining such identities.

Remind that the *binomial transform* of a sequence $a_1, a_2, \ldots$, is a new sequence $b_1, b_2, \ldots$, such that for every positive integer $k$,

$$b_k = \sum_{j=1}^{k}(-1)^{k-j}\binom{k}{j} a_j, \text{ with inversion } a_k = \sum_{j=1}^{k}\binom{k}{j} b_j \qquad (2)$$

([8, (5.48), p. 192], [9], [10]). Equation (1) shows that the sequences $k!S(\alpha,k)$ and $k^\alpha$ are related by the binomial transform. The inversion formula then yields

$$k^\alpha = \sum_{j=1}^{k} \binom{k}{j} j!\, S(\alpha,j),\qquad(3)$$

for any positive integer $k$.

## 2. The Identities

We start with a simple lemma

**Lemma.** Let $c_1, c_2, \ldots$, be a sequence of complex numbers. Then for every positive integer $m$ we have

$$\sum_{k=1}^{m} k^\alpha c_k = \sum_{j=1}^{m} j!\, S(\alpha,j) \sum_{k=j}^{m} \binom{k}{j} c_k . \qquad(4)$$

For the proof we just need to use (3) for $k^\alpha$ and then change the order of summation on the right hand side

$$\sum_{k=1}^{m} k^\alpha c_k = \sum_{k=1}^{m} c_k \sum_{j=1}^{k} \binom{k}{j} j!\, S(\alpha,j) = \sum_{j=1}^{m} j!\, S(\alpha,j) \sum_{k=j}^{m} \binom{k}{j} c_k . \qquad(5)$$

This lemma helps to generate power sum identities by using various upper summation identities. We present here five examples arranged in four propositions.

**Proposition 1.** For every positive integer $m$ and every two complex numbers $\alpha \neq 0, x$,

$$\sum_{k=1}^{m} k^\alpha x^k = \sum_{j=1}^{m} j!\, S(\alpha,j)\, \sigma(x,m,j), \qquad(6)$$

where $\sigma(x,m,j)$ is the (upper summation) polynomial

$$\sigma(x,m,j) = \sum_{k=j}^{m} \binom{k}{j} x^k = x^j \sum_{r=0}^{m-j} \binom{r+j}{j} x^r . \qquad(7)$$

In particular, when $x = 1$ one has

$$\sum_{k=1}^{m} k^\alpha = \sum_{j=1}^{m} \binom{m+1}{j+1} j!\, S(\alpha,j) . \qquad(8)$$

*Proof.* We use the Lemma with $c_k = x^k$. When $x = 1$ we use the upper summation identity

$$\sum_{k=j}^{m} \binom{k}{j} = \binom{m+1}{j+1} \tag{9}$$

(see, for instance, [7, 1.52] or [8, p.174]) . Thus (6) turns into (8).

**Remark.** Identity (8) was proved in [2] in the equivalent form

$$\sum_{k=1}^{m} k^\alpha = \sum_{j=1}^{m} \binom{m}{j} (j-1)!\, S(\alpha+1, j) \tag{10}$$

by induction. The equivalence follows from the properties

$$S(\alpha+1, k) = k S(\alpha, k) + S(\alpha, k-1) \tag{11}$$

(see [2, (1.16)]), and the well known binomial identity [8, p.174],

$$\binom{m}{k} + \binom{m}{k-1} = \binom{m+1}{k}. \tag{12}$$

**Remark.** With complex powers $\alpha \neq 0$ as in (6) we have the flexibility to write

$$\sum_{k=1}^{m} \frac{x^k}{k^\alpha} = \sum_{j=1}^{m} j!\, S(-\alpha, j)\, \sigma(x, m, j). \tag{13}$$

When $\alpha = n$ is a positive integer, identity (8) (or (10), to that matter) is well known and has a long history. In the early 18-th century, Bernoulli evaluated $\sum_{k=1}^{m} k^n$ in terms of the numbers known today as Bernoulli numbers. Continuing Bernoulli's work, Leonhard Euler [4, paragraphs 173, 176] evaluated sums of the form $\sum_{k=1}^{m} k^n x^k$, essentially by applying $n$ times the operator $x \dfrac{d}{dx}$ to the identity

$$\sum_{k=1}^{m} x^k = \frac{1}{1-x} - \frac{x^{m+1}}{1-x} \tag{14}$$

$(x \neq 1)$. This led him to the discovery of a special sequence of polynomials $A_k(x)$ called today

Eulerian polynomials ([1], [3], [6]). In terms of these polynomials one has

$$(x\frac{d}{dx})^n \frac{1}{1-x} = \frac{A_n(x)}{(1-x)^{n+1}}, \quad n = 0, 1, \ldots, \tag{15}$$

and therefore, with some help from the Leibniz rule

$$\sum_{k=1}^{m} k^n x^k = \frac{A_n(x)}{(1-x)^{n+1}} - x^{m+1} \sum_{k=0}^{n} \binom{n}{k} \frac{(m+1)^{n-k} A_k(x)}{(1-x)^{k+1}}. \tag{16}$$

This identity, however, cannot be extended to complex powers $n \to \alpha \in \mathbb{C}$ for obvious reasons.

The next identity can be viewed as the binomial transform of the sequence $k^\alpha x^k$ extending equation (1).

**Proposition 2.** For every positive integer $m$ and every two complex numbers $\alpha \neq 0, x$,

$$\sum_{k=1}^{m} \binom{m}{k} k^\alpha x^k = \sum_{j=1}^{m} \binom{m}{j} j! \, S(\alpha, j) x^j (1+x)^{m-j}. \tag{17}$$

*Proof.* We apply the Lemma with $c_k = \binom{m}{k} x^k$. The result then follows from the interesting identity

$$\sum_{k=j}^{m} \binom{m}{k} \binom{k}{j} x^k = \binom{m}{j} x^j (1+x)^{m-j}, \tag{18}$$

which is listed as number 3.118 on p. 36 in [7]. To prove this identity one can start by reducing both sides by $x^j$ and then expanding $(1+x)^{m-j}$.

Note that when $x = -1$, (17) turns into (1).

**Remark.** Identity (17) for positive integers $\alpha = r$ can also be found in the treasure chest [7]. It is listed there (as number 1.126 on p.16) in the form

$$\sum_{k=0}^{n} \binom{n}{k} k^r x^k = (1+x)^n \sum_{j=0}^{r} (-1)^j \binom{n}{j} \frac{x^j}{(1+x)^j} \sum_{k=0}^{j} (-1)^k \binom{j}{k} k^r. \tag{19}$$

Note that in (19) the number $r$ *has* to be a positive integer, because it stands for the upper limit of the first sum on the RHS. For the case $x = 1$, (19) was recently rediscovered by Spivey [10].

The next identity involves the unsigned Stirling numbers of the first kind $\begin{bmatrix} n \\ k \end{bmatrix}$ (see [8]).

**Proposition 3**. For every positive integer $m$ and every complex $\alpha \neq 0$ we have

$$\sum_{k=1}^{m} \begin{bmatrix} m \\ k \end{bmatrix} k^\alpha = \sum_{j=1}^{m} j! S(\alpha, j) \begin{bmatrix} m+1 \\ j+1 \end{bmatrix}. \qquad (20)$$

The proof uses the lemma with $c_k = \begin{bmatrix} m \\ k \end{bmatrix}$ and also the upper summation identity ([8, (6.16), p. 265])

$$\sum_{k=j}^{m} \binom{k}{j} \begin{bmatrix} m \\ k \end{bmatrix} = \begin{bmatrix} m+1 \\ j+1 \end{bmatrix}. \qquad (21)$$

We finish this note with two identities involving the harmonic numbers

$$H_k = 1 + \frac{1}{2} + \ldots + \frac{1}{k}, \quad (k = 1, 2, \ldots). \qquad (22)$$

**Proposition 4**. For every positive integer $m$ and every complex power $\alpha \neq 0$,

$$\sum_{k=1}^{m} H_k k^\alpha = \sum_{j=1}^{m} j! S(\alpha, j) \binom{m+1}{j+1} (H_{m+1} - \frac{1}{j+1}), \qquad (23)$$

$$\sum_{k=1}^{m} \frac{k^\alpha}{m-k+1} = \sum_{j=1}^{m} j! S(\alpha, j) \binom{m+1}{j} (H_{m+1} - H_j). \qquad (24)$$

*Proof*. Follows from the lemma with $c_k = H_k$ and $c_k = \frac{1}{m+k-1}$ correspondingly and also from the two upper summation identities - see [8, (6.70), p. 280 and p. 354],

$$\sum_{k=j}^{m} \binom{k}{j} H_k = \binom{m+1}{j+1} (H_{m+1} - \frac{1}{j+1}) \qquad (25)$$

$$\sum_{k=j}^{m} \binom{k}{j} \frac{1}{m-k+1} = \binom{m+1}{j} (H_{m+1} - H_j). \qquad (26)$$